\documentclass[a4paper,11pt]{article}

\usepackage{subfig}
\usepackage{pgfplots}
\usepackage{pgfplotstable}
\usepackage{mathtools}
\usepackage{multicol}
\usepackage{comment}
\usepackage{booktabs}
\pgfplotsset{compat=1.5}
\usepackage{amssymb}
\usepackage{url}
\usepackage{bm}

\usepackage{pdfsync}
\usepackage{float}
\usepackage{tabularx}
\usepackage{enumerate}
\usepackage{array}
\usepackage{xspace}

\usepackage{authblk}

\newcommand{\x}{\ensuremath{\boldsymbol{x}}}
\newcommand{\mupar}{\ensuremath{\boldsymbol{\mu}}}

\begin{document}

\title{A complete data-driven framework for the efficient solution of
  parametric shape design and optimisation in naval engineering
  problems}

\author[]{Nicola~Demo\footnote{nicola.demo@sissa.it}}
\author[]{Marco~Tezzele\footnote{marco.tezzele@sissa.it}}
\author[]{Andrea~Mola\footnote{andrea.mola@sissa.it}}
\author[]{Gianluigi~Rozza\footnote{gianluigi.rozza@sissa.it}}

\affil[]{Mathematics Area, mathLab, SISSA, International School of Advanced Studies, via Bonomea 265, I-34136 Trieste, Italy}

\maketitle

\begin{abstract}
In the reduced order modeling (ROM) framework, the solution of a
parametric partial differential equation is approximated by combining the
high-fidelity solutions of the problem at hand for several properly chosen
configurations.  Examples of the ROM application, in the naval field, can be
found in~\cite{tezzele2018ecmi,salmoiraghi2016advances}. Mandatory ingredient
for the ROM methods is the relation between the high-fidelity solutions and the
parameters.  Dealing with geometrical parameters, especially in the industrial
context, this relation may be unknown and not trivial (simulations over hand
morphed geometries) or very complex (high number of parameters or many nested
morphing techniques). To overcome these scenarios, we propose in this
contribution an efficient and complete data-driven framework involving ROM
techniques for shape design and optimization, extending the pipeline presented
in~\cite{demo2018shape}.  By applying the singular value decomposition (SVD) to
the points coordinates defining the hull geometry --- assuming the topology is
inaltered by the deformation ---, we are able to compute the optimal space
which the deformed geometries belong to, hence using the modal coefficients as
the new parameters we can reconstruct the parametric formulation of the domain.
Finally the output of interest is approximated using the proper orthogonal
decomposition with interpolation technique. To conclude, we apply this framework to a
naval shape design problem where the bulbous bow is morphed to reduce the total
resistance of the ship advancing in calm water.
\end{abstract}



\section{Introduction}
\label{sec:intro}

The reduced basis method (RBM)~\cite{HesthavenRozzaStamm2015,rozza2007reduced}
is a well-spread technique for reduced order modeling, both in academia
and in industry~\cite{salmoiraghi2016advances,rozza2018advances,tezzele2018ecmi,morhandbook2019},
and consists in two phases: an offline phase that can be carried out
on high performance computing facilities, and an online one that
exploits the reduced dimensionality of the system to perform the
parametric computation on portable devices. In the 
offline stage the reduced order space is created from full order complex
simulations computed for certain values of the parameters. The selection of the
reduced basis functions that span this new reduced space can be carried out by
different techniques. In this work we employ the proper orthogonal
decomposition (POD)~\cite{quarteroni2014,benner2017}, which is based on the
singular value decomposition (SVD), on the set of high-fidelity snapshots.
After the creation of such space, in the online phase a new parametric solution
is calculated as a linear combination of the precomputed reduced basis
functions. The creation of a reduced order model is crucial in the
shape optimisation context where the optimiser needs to compute
several high-fidelity simulations. 

Novelty of this work is the creation of a reduced order space
containing the manifold of admissible shapes by applying POD over the sampled
geometries, in order to reduce the parameter space dimension and
to enhance the order reduction of the output fields. To generate the
original design space we employ the free form deformation (FFD)
method, a well-known shape parametrisation technique. Another approach
for reduced order modeling enhanced by parameter space reduction technique can be found
in~\cite{tezzele2018combined} where they propose a coupling between
POD-Galerkin methods and active subspaces. 
After the creation of the reduced space for the admissible shapes, we
can exploit this new parametric formulation for the construction of
the reduced space for the output fields, using the non-intrusive
technique called POD with interpolation
(PODI)~\cite{bui2003proper,Carlberg2010,everson1995karhunen} for the 
online computation of the coefficients of the linear combination.
We would like to cite~\cite{meng2018nonlinear} where they present the concept of a
shape manifold representing all the admissible shapes, independently
of the original design parameters, and thus exploiting the intrinsic
dimensionality of the problem.

This work is organised as follows: after the presentation of the
general setting of the problem, there is a brief overview of the FFD
method, then we illustrate how the parameter space reduction is
performed, and we present PODI for the reduction of the high-fidelity
snapshots. Finally the numerical results are presented with the
conclusions and some perspective.

\section{The problem}
\label{sec:problem}

\begin{figure}[h!]
\centering
\includegraphics[trim=0 0 0 0, width=.75\textwidth, height=3cm, keepaspectratio]{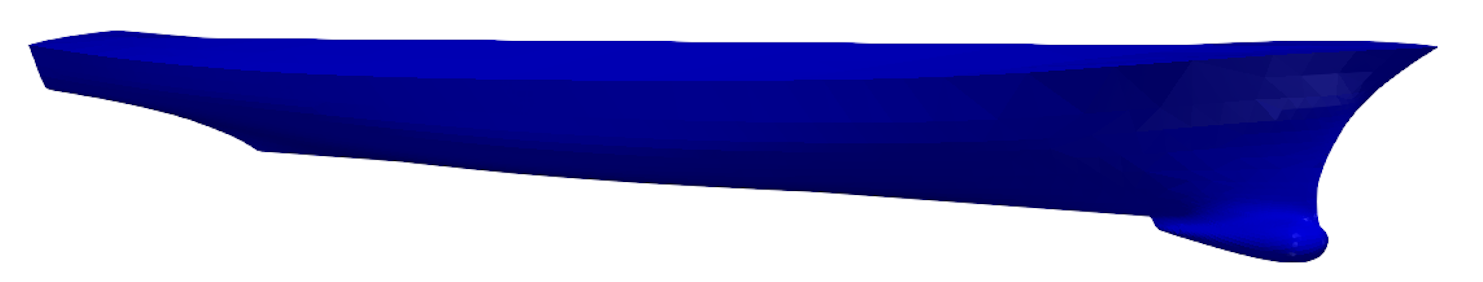}\hfill
\includegraphics[trim=0 0 0 0, width=.24\textwidth, height=3cm, keepaspectratio]{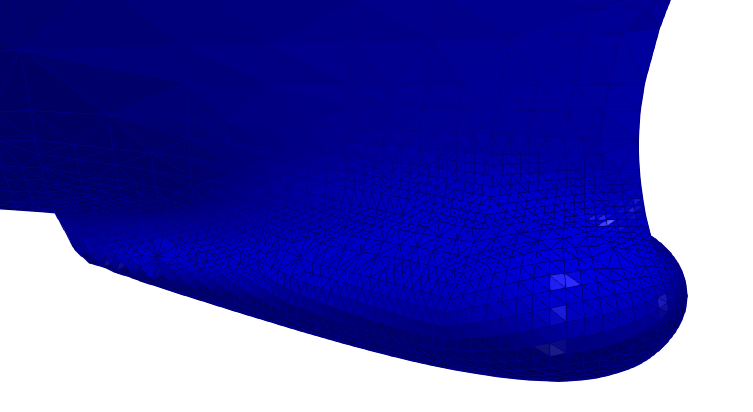}
\caption{Complete hull domain representing the DTMB 5415 and, on the right, a zoom on the bulbous
  bow.}
\label{fig:domain}
\end{figure}

Let $\Omega \subset \mathbb{R}^3$ be the reference hull domain.
We define a parametric shape morphing function $\mathcal{M}$ as follows
\begin{equation}
  \label{eq:general_morphing}
  \mathcal{M}(\boldsymbol{x}; \mupar): \mathbb{R}^3 \to \mathbb{R}^3,
\end{equation}
which maps $\Omega$ into the deformed domain $\Omega(\mupar)$ as 
$\Omega(\mupar) = \mathcal{M}(\Omega; \mupar)$, 
where $\mupar \in \mathbb{D} \subset \mathbb{R}^5$ represents the
vector of the geometrical parameters. $\mathbb{D}$ will be properly defined in
Section~\ref{sec:parameters}. Such map $\mathcal{M}$ can represent
many different morphing techniques (not necessarily affine) such as free form
deformation (FFD)~\cite{sederbergparry1986}, radial basis functions (RBF)
interpolation~\cite{buhmann2003radial,morris2008cfd,manzoni2012model},
and the inverse distance weighting (IDW)
interpolation~\cite{shepard1968,witteveenbijl2009,forti2014efficient,BallarinDAmarioPerottoRozza2018},
for instance. In this work we use the FFD, presented in
Section~\ref{sec:ffd}, to morph a bulbous bow of a benchmark hull. We
chose the DTMB~5415 hull thanks to the vast amount of experimental data available in
the literature, see for example~\cite{olivieri2001towing}. In
Figure~\ref{fig:domain} the domain $\Omega$ and a particular of the
bulbous bow we are going to parametrize and deform.

The pipeline is the following: using geometrical FFD parameters we
generate several deformed hulls; then we apply the POD on the coordinates
of the points describing the deformed geometries, and the new
parameters will be the POD coefficients of the selected modes; after
checking for possible linear dependencies between these coefficients,
we sample the reduced parameter space producing new deformed hulls
upon which we are going to actually perform CFD simulations.
Regarding the full order model, we use the Reynolds-averaged Navier-Stokes
equations to describe the incompressible and turbulent flow around the ship.
The Froude number has been set to 0.2 and  we chose the $k$--$\omega$ SST
model for the turbulence since it is one of the most popular benchmark for
hydrodynamic analysis for industrial naval problems.
In this way we reduce the parameter space, staying on the manifold of the
admissible shapes, and reducing the burden of the output reduced space
construction through PODI.

\section{Free form deformation of the bulbous bow}
\label{sec:ffd}

Here we are going to properly define the deformation map $\mathcal{M}$
introduced in Eq.~\eqref{eq:general_morphing}, which we employed for
this work, and that corresponds to the free form deformation (FFD)
technique.
The original formulation of the FFD can be found
in~\cite{sederbergparry1986}, for more
recent works in the context of reduced basis methods for shape
optimization we
cite~\cite{LassilaRozza2010,rozza2013free,sieger2015shape}. It has
also been applied to naval engineering problems
in~\cite{demo2018shape,demo2018isope,tezzele2018dimension}, 
while for an automotive case see~\cite{salmoiraghi2018}.

\begin{figure}[h!]
\centering
\includegraphics[trim=0 0 0 0, width=.60\textwidth]{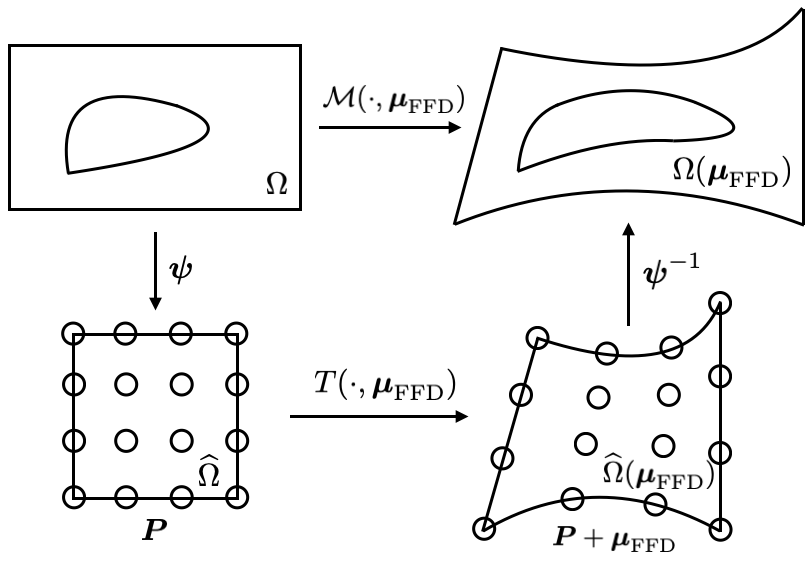}
\caption{Sketch of the FFD map $\mathcal{M}$ composition. The domain
  is mapped to a reference configuration, then the lattice of FFD control
  points induce the body deformation, and finally the morphed object
  is mapped back to the physical space.}
\label{fig:ffd_scheme}
\end{figure}

The FFD map is the composition of three maps described in the
following, while for a visual representation we refer to Figure~\ref{fig:ffd_scheme}:
\begin{itemize}
\item the function $\boldsymbol{\psi}$ maps the physical domain to the
  reference one where we construct the reference lattice of points,
  denoted with $\boldsymbol{P}$ around the object to be morphed;
\item the function $T$ performs the actual deformation since it
  applies the displacements defined by $\mupar_{\text{FFD}}$ to the
  lattice $\boldsymbol{P}$. It uses the B-splines or Bernstein
  polynomials tensor product to morph all the points inside the
  lattice of control points;
\item finally we need to map back the deformed domain to the physical
  configuration through the map $\boldsymbol{\psi}^{-1}$.
\end{itemize}

Se we can define the FFD map $\mathcal{M}$ through the composition of
the three maps presented above as
\begin{equation}
\mathcal{M} (\x, \mupar_{\text{FFD}}) := (\boldsymbol{\psi}^{-1} \circ T
\circ \boldsymbol{\psi}) (\x, \mupar_{\text{FFD}})
\quad \forall \x \in \Omega.
\end{equation}

\begin{figure}[h!]
\centering
\includegraphics[trim=0 0 0 0, width=.60\textwidth]{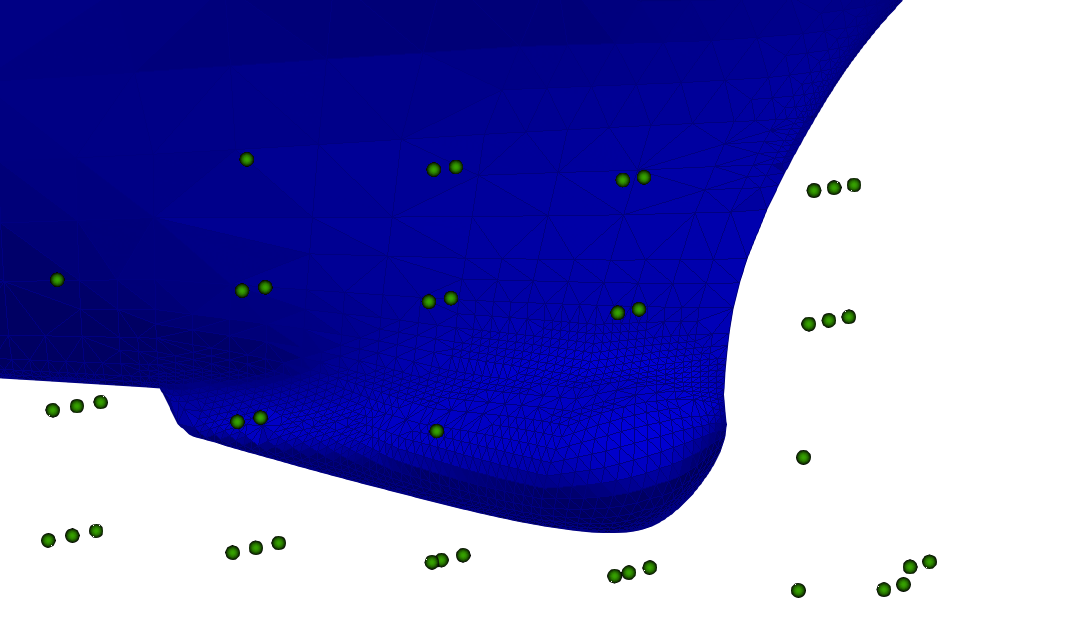}
\caption{Example of FFD parametrisation and morphing of the DTMB~5415
  hull. In green the lattice of control points that define the actual deformation.}
\label{fig:ffd_lattice}
\end{figure}

In Figure~\ref{fig:ffd_lattice} it is possible to see the actual
lattice of points we used, in green, for a particular choice of the
FFD parameters. For an actual implementation of this method in Python,
along with other possibile deformation methods, we refer to the open
source package called PyGeM - Python Geometrical Morphing~\cite{pygem}.

\section{Reduction of the parameter space through POD of the mesh coordinates}
\label{sec:parameters}

In order to reduce the parameter space dimension we apply the
POD on a set of snapshots that depends on the FFD
parameters. Each snapshot is the collection of all the coordinates
of the points defining the \emph{stl} file geometry. Since the
generation of these snapshots does not depend on complex simulations
but only on the particular FFD deformation, we are able to create a
dataset with as many entries as we want.
So we create a database of $N_{\text{train}} = 1500$ geometrical parameters
$\mupar_{\text{FFD}} \in \mathbb{D} := [-0.3, 0.3]^5$ sampled 
with a uniform distribution. Moreover we create the corresponding
database of mesh coordinates $\mathbf{u}$ corresponding to these parameters, that
is $\boldsymbol{\Theta} = [\mathbf{u} (\mupar_{\text{FFD}, \, 1}) | \dots |
\mathbf{u} (\mupar_{\text{FFD}, \, N_{\text{train}}})]$. Then we perform the singular
value decomposition (SVD) on $\boldsymbol{\Theta}$ in order to extract
the matrix of POD modes:
\begin{equation}
\boldsymbol{\Theta} = \boldsymbol{\Psi} \boldsymbol{\Sigma} \boldsymbol{\Phi}^T,
\end{equation}  
where with $\boldsymbol{\Psi}$ and $\boldsymbol{\Phi}$ we denote the
left and right singular vectors matrices of $\boldsymbol{\Theta}$
respectively, and with $\boldsymbol{\Sigma}$ the diagonal matrix
containing the singular values in decreasing order. The columns of
$\boldsymbol{\Psi}$, denoted with $\psi_i$, are the so-called POD modes.
We can thus express the approximated reduced mesh with the first $N$
modes as
\begin{equation}
\mathbf{u}^N = \sum_{i=1}^N \alpha_i \psi_i ,
\label{eq:lcombi}
\end{equation}
where $\alpha_i$ are the so called POD coefficients. To compute them
in matrix form we just use the database we created as follows
\begin{equation}
\boldsymbol{\alpha} = \boldsymbol{\Psi}^T \boldsymbol{\Theta},
\end{equation}
and then we truncate to the first $N$ modes and
coefficients.

\begin{figure}[h!]
\centering
\includegraphics[trim=20 0 70 20, width=.49\textwidth]{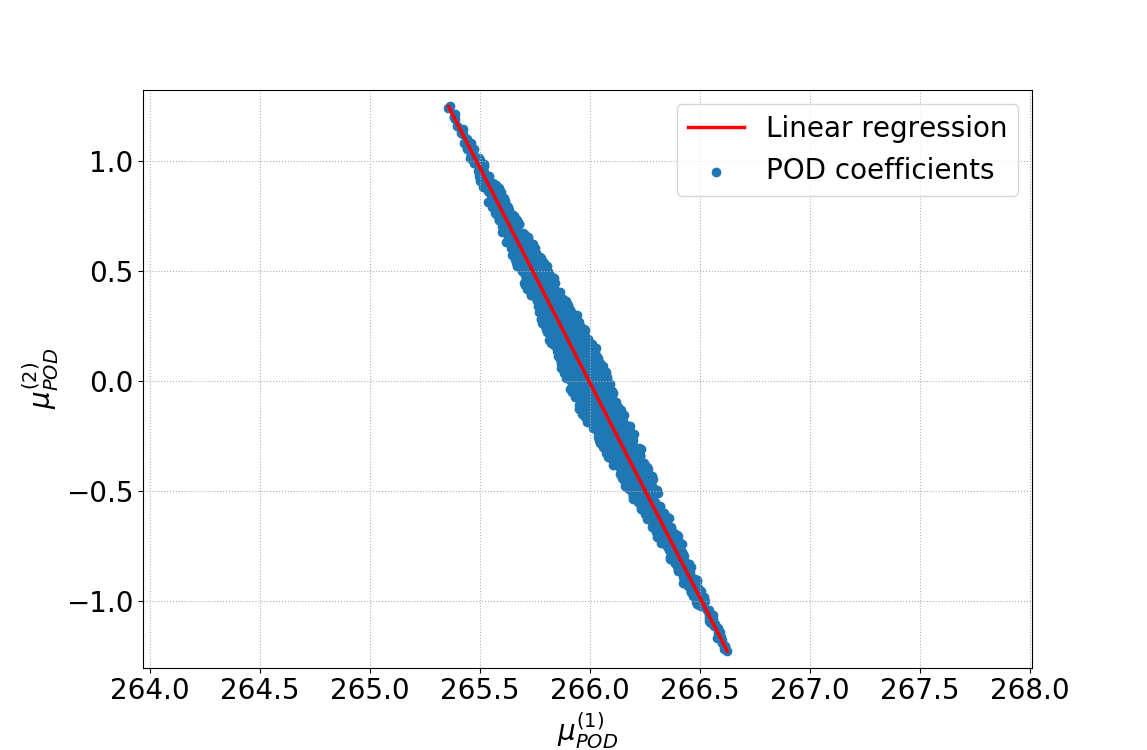}\hfill
\includegraphics[trim=20 0 70 20, width=.49\textwidth]{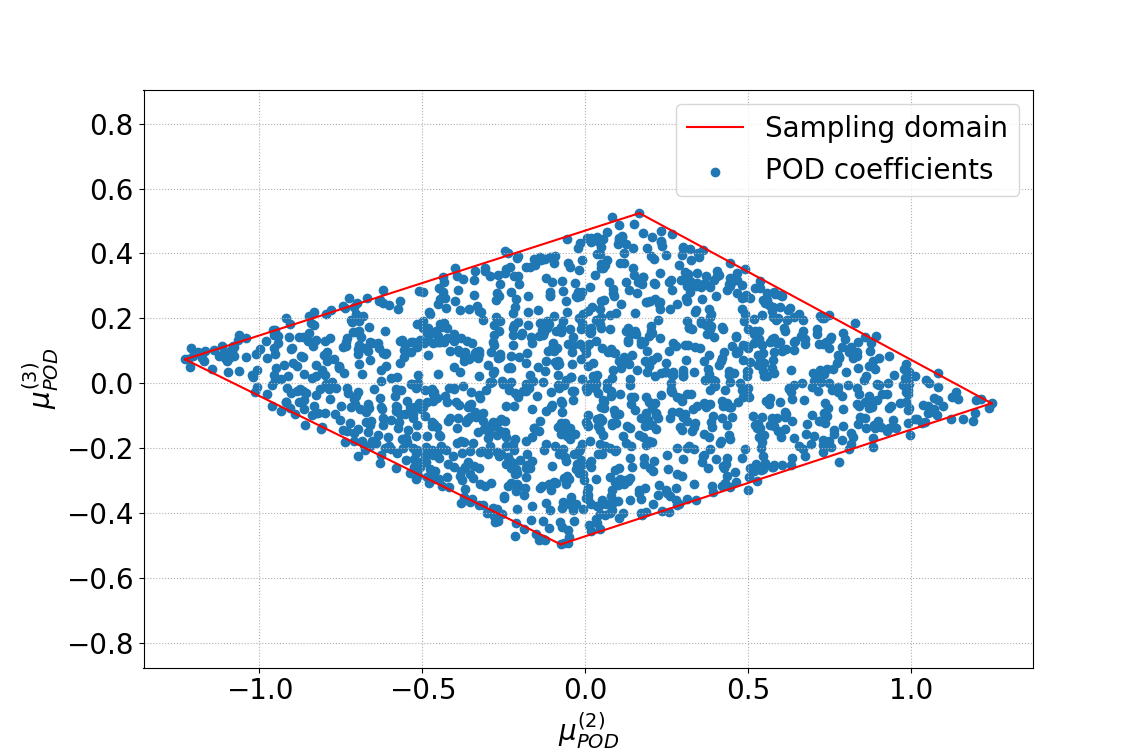}
\caption{POD coefficients dependance. On the left we have
  $\mupar_{\text{POD}}^{(2)}$ with respect to
  $\mupar_{\text{POD}}^{(1)}$ in blue, and in red the linear
  regression to approximate one as a function of the other. On the
  right $\mupar_{\text{POD}}^{(3)}$ against
  $\mupar_{\text{POD}}^{(2)}$ in blue, and in red the
  boundaries defining the quadrilateral in which the sampling is
  performed.}
\label{fig:pod_coeffs}
\end{figure}

After the selection of the number of POD modes required to have an
accurate approximation of each geometry, we end up with the first
reduction of the parameter space, that is with 3 POD coefficients
$\mupar_{\text{POD}} := \boldsymbol{\alpha} \in \mathbb{R}^3$, we are
able to represent all the possible deformations for
$\mupar_{\text{FFD}} \in \mathbb{D}$. So we can express every geometry
with 3 modes, but the coefficients can still be linearly dependent. We
can investigate this dependance by plotting every component
$\mupar_{\text{POD}}^{(i)}$ against each  
other. As we can see from the plot on the left in
Figure~\ref{fig:pod_coeffs}, we can approximate
$\mupar_{\text{POD}}^{(2)}$ with a linear regression given
$\mupar_{\text{POD}}^{(1)}$. For what concerns
$\mupar_{\text{POD}}^{(3)}$, we can constraint it to be inside the
quadrilateral in Figure~\ref{fig:pod_coeffs}, on the right. So we are able to
express every possible geometry described with the original 5 FFD
parameters with only 2 new independent parameters. We can thus sample
the full parameter space using a new reduced space, preserving the
geometrical variability, and reducing the construction cost of the
reduced output field space. This, as we are going to present, results
in a faster optimization procedure.

\section{Non-intrusive reduced order modeling by means of PODI}
\label{sec:podi}

Proper orthogonal decomposition with interpolation is a non-intrusive data-driven
method for reduced order modeling allowing an efficient approximation of the
solution of parametric partial differential equations.  As well as for the
geometries, we collect in a database the high-fidelity solutions of several CFD
simulations corresponding to different configurations, then we apply the POD
algorithm to the solutions matrix --- the matrix whose columns are the
solutions --- in order to extract the POD modes that span the optimal space
which the solutions belong to. Thus the solutions can be projected onto the
reduced space: we represent the high-fidelity solutions as linear combination
of the POD modes. Similarly to Eq.~\eqref{eq:lcombi}, the modal coefficients
of the $i$-th solution
$\mathbf{x}^{\text{PODI}}_i$ --- also called the reduced solution --- are
obtained as:
\begin{equation}
\mathbf{x}^{\text{PODI}}_i = \mathbf{U}^T\mathbf{x}_i\qquad \forall i
\in \{1, \dotsc, M\}
\end{equation}
where $\mathbf{U}$ refers to the POD modes and $M$ is the number of
high-fidelity solutions. We call $N$ the number of POD modes and $\mathcal{N}$
the dimension of high-fidelity solutions then $\mathbf{x}^{\text{PODI}}_i \in V^N$ and
$\mathbf{x}_i \in V^{\mathcal{N}}$. Since in complex problems we have an high
number of degrees of freedom, typically we have $N \ll \mathcal{N}$. The low-rank
representation of the solutions allows to easily interpolate them, exploting
the relation between the reduced solutions and the input
parameters: in this way, we can compute the modal
coefficients for any new parametric point and project the reduced solution onto
the high dimensional space for a real-time approximation of the truth solution.
This technique is defined non-intrusive, since it relies only on the solutions,
without requiring information about the physical system and the equations
describing it. For this reason it is particularly suited for industrial
problem, thanks to its capability to be coupled also with commercial solvers. The
downside is the error introduced by the interpolation, depending by the method
itself, and the requirement of solutions with the same dimensionality, that can
be a problem if the computational grid is built from scratch for any new
configuration. Possible solutions are the projection of the solution on a
reference mesh~\cite{demo2018shape}, or to deform the grid using the laplacian
diffusion~\cite{stabile2019}. Moreover, we cite~\cite{garotta2018quiet,
salmoiraghi2018} for other examples of PODI applications.  For this work, we
employed the open source Python package EZyRB~\cite{demo2018ezyrb} as
software to perform the data-driven model order reduction.

\section{Numerical results}
\label{sec:results}

In this section we present the results for the application of the
complete pipeline to the problem presented in
Section~\ref{sec:problem}.

\begin{figure}[h!]
\centering
\includegraphics[trim=0 0 0 0, width=.8\textwidth]{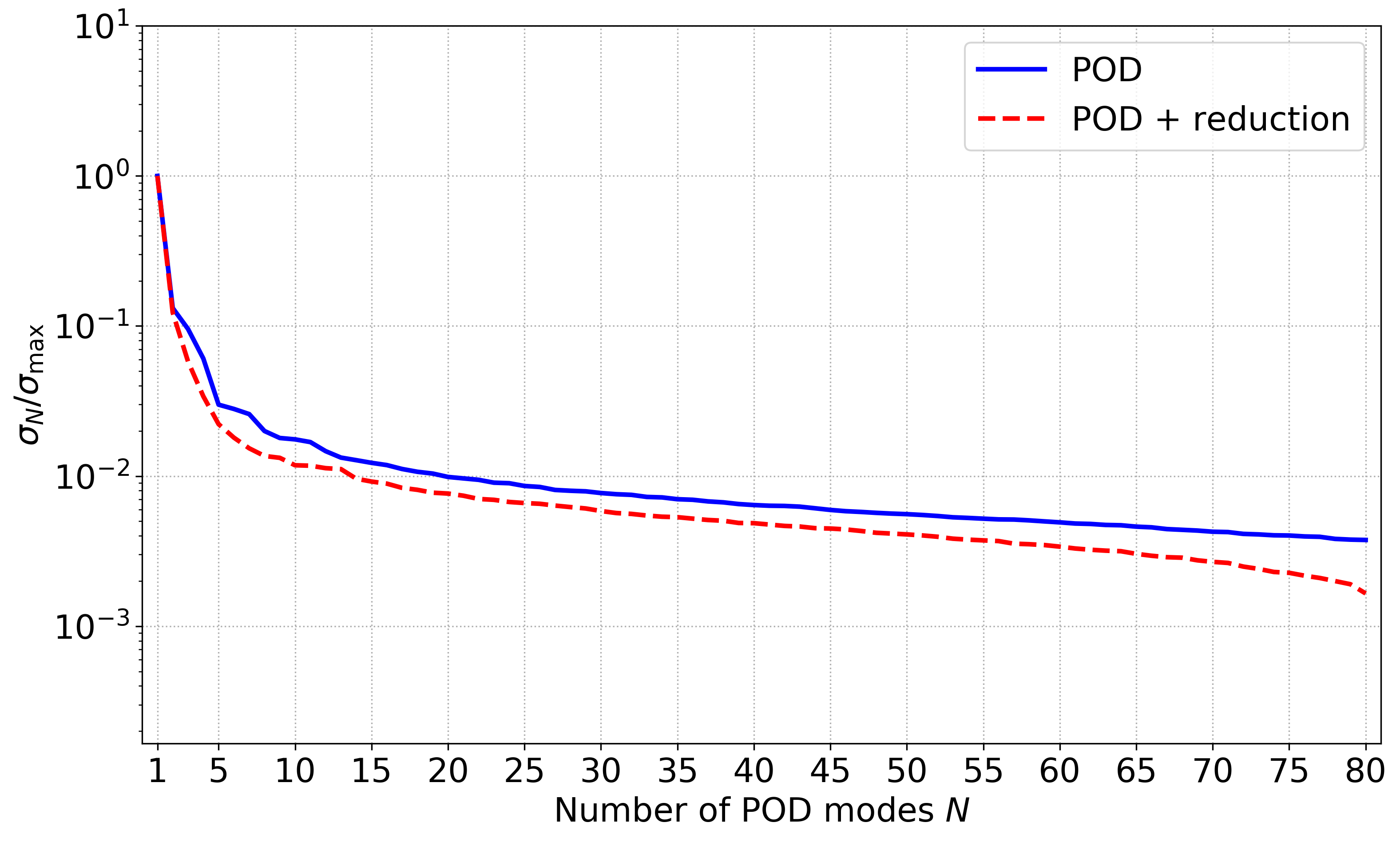}
\caption{POD singular values decay as a function of the number of
  modes. The blue line corresponds to the original sampling in the
  full parameter space, while the red dotted line, which identifies
  the POD+reduction approach, corresponds to the sampling in the new
  reduced parameter space.}
\label{fig:svd}
\end{figure}

First, we sample the full parameter space $\mathbb{D}$ extracting
$N_{\text{POD}} = 100$ parameters to construct the reduce order model
without any further reduction, and we identify this approach with the
subscript ``POD''. Then, as explained in Section~\ref{sec:parameters},
we compute the shape manifold with 1500 different deformations, and we
extract the new coefficients describing the new reduced parameter
space. We sample this 2-dimensional space uniformly and we collect
$N_{\text{POD+reduction}} = 80$ solution snapshots. We can compare the
decay of the singular values of the snapshots matrix for the two
approaches. In Figure~\ref{fig:svd} we can note how the proposed
computational pipeline results in a faster decay and thus in a better
approximation for a given number of POD modes.

We underline that, despite the gain is not so big, the results do not
involve further high-fidelity simulations. We only collected several
different deformations at a negligible computational cost with respect
to a single full order CFD simulation. Moreover the construction of
the interpolator takes a huge advantage of the reduced parameter space
since it counters the curse of dimensionality.

We can conclude that the proposed preprocessing step has sever
benefits in terms of accuracy of the reduced order model at a small
cost from a computational point of view.

\section{Conclusions and perspectives}
\label{sec:the_end}

In this work we presented a complete data-driven numerical pipeline
for shape optimization in naval engineering problems. The object was
to find the optimal bulbous bow to minimize the total drag resistance
of a hull advancing in calm water. First we parametrized and morphed
the bulbous bow through the free form deformation method. Then we
reduced the parameter space dimension approximating the shape manifold
with the use of proper orthogonal decomposition and the investigation
on linear dependance of the POD coefficients. We create the reduced
order model sampling only the reduced two dimensional parameter space
and with POD with interpolation we can compute in real time the outputs
of interest for untried new parameters. Thus the optimizer can query
the surrogate model and find the optimal shape.

\section*{Acknowledgements}
This work was partially performed in the context of the project SOPHYA -
``Seakeeping Of Planing Hull YAchts'', supported by Regione
FVG, POR-FESR 2014-2020, Piano Operativo Regionale Fondo Europeo per
lo Sviluppo Regionale, and partially supported by European Union Funding for
Research and Innovation --- Horizon 2020 Program --- in the framework
of European Research Council Executive Agency: H2020 ERC CoG 2015
AROMA-CFD project 681447 ``Advanced Reduced Order Methods with
Applications in Computational Fluid Dynamics'' P.I. Gianluigi
Rozza.


\end{document}